\documentclass[12pt]{article}%
\usepackage{graphicx}
\usepackage[intlimits]{amsmath}
\usepackage{latexsym}
\usepackage{amsfonts}
\usepackage{amssymb}%
\makeatletter
\newfont{\footsc}{cmcsc10 at 8truept}
\newfont{\footbf}{cmbx10 at 8truept}
\newfont{\footrm}{cmr10 at 10truept}
\newtheorem{theorem}{Theorem}

\newtheorem{proposition}[theorem]{Proposition}

\newenvironment{proof}[1][Proof]{\noindent{\textbf {#1}  }}  {\hfill$\Box$\bigskip}

\begin{document}

\title{Revisiting Schur's bound on the largest singular value}
\author{Vladimir Nikiforov\\Department of Mathematical Sciences, University of Memphis, \\Memphis TN 38152, USA, email: \textit{vnikifrv@memphis.edu}}
\maketitle

\begin{abstract}
We give upper and lower bounds on the largest singular value of a matrix using
analogues to walks in graphs. For nonnegative matrices these bounds are
asymptotically tight.

In particular, the following result improves a bound due to Schur. If
$A=\left(  a_{ij}\right)  $ is an $m\times n$ complex matrix, its largest
singular value satisfies%
\[
\sigma^{2}\left(  A\right)  \leq\max_{i\in\left[  m\right]  }%
{\textstyle\sum\limits_{j\in\left[  n\right]  }}
\left\vert a_{ij}\right\vert c_{j}\leq\max_{a_{ij}\neq0}r_{i}c_{j},
\]
where $r_{i}=\sum_{k\in\left[  n\right]  }\left\vert a_{ik}\right\vert ,$
$c_{j}=\sum_{k\in\left[  m\right]  }\left\vert a_{kj}\right\vert .\medskip$



\end{abstract}

\section{Introduction}

How large the largest singular value $\sigma\left(  A\right)  $ of an $m\times
n$ matrix $A=\left(  a_{ij}\right)  $ can be? In 1911 Schur \cite{Sch11}, p.
6, gave the bound
\begin{equation}
\sigma^{2}\left(  A\right)  \leq\max_{i\in\left[  m\right]  ,\text{ }%
j\in\left[  n\right]  }r_{i}c_{j}, \label{Schin}%
\end{equation}
where $r_{i}=\sum_{k\in\left[  n\right]  }\left\vert a_{ik}\right\vert ,$
$c_{j}=\sum_{k\in\left[  m\right]  }\left\vert a_{kj}\right\vert .$

The aim of this note to strengthen this bound and give similar lower bounds on
$\sigma\left(  A\right)  .$ In particular, our results imply that if $A$ is
nonzero, then
\begin{equation}
\sigma^{2}\left(  A\right)  \leq\max_{i\in\left[  m\right]  }%
{\textstyle\sum\limits_{j\in\left[  n\right]  }}
\left\vert a_{ij}\right\vert c_{j}\leq\max_{a_{ij}\neq0}r_{i}c_{j}.
\label{ourin}%
\end{equation}
Note that sometimes (\ref{ourin}) is much stronger than (\ref{Schin}). Indeed,
letting $A$ be the adjacency matrix of the star $K_{1,n},$ inequality
(\ref{Schin}) gives $\sigma^{2}\left(  A\right)  \leq n^{2},$ while
(\ref{ourin}) gives $\sigma^{2}\left(  A\right)  \leq n,$ which is best
possible, in view of $\sigma^{2}\left(  A\right)  =n$.

For basic notation and definitions see \cite{HoJo88}. In particular,
$\mathbf{j}_{m}$ denotes the vector of $m$ ones.

Given an $m\times n$ matrix $A=\left(  a_{ij}\right)  $, for all $r\geq0$ and
$i,j\in\left[  m\right]  ,$ let $w_{A}^{r}\left(  i,j\right)  $ be the
$\left(  i,j\right)  $th entry of $\left(  AA^{\ast}\right)  ^{r}.$ Set
$w_{A}^{r}\left(  k\right)  =%
{\textstyle\sum_{i\in\left[  m\right]  }}
w_{A}^{r}\left(  k,i\right)  $ and $w_{A}^{r}=%
{\textstyle\sum_{i\in\left[  m\right]  }}
w_{A}^{r}\left(  i\right)  .$

Note that if $A$ is the adjacency matrix of a graph, then $w_{A}^{r}$ is the
number of walks on $2r+1$ vertices$.$

The following theorem generalizes inequality (\ref{ourin}) and thus,
inequality (\ref{Schin}).

\begin{theorem}
\label{th1}For every nonzero $m\times n$ matrix $A=\left(  a_{ij}\right)  $
and all $r\geq0,$ $p\geq1,$
\[
\sigma^{2p}\left(  A\right)  \leq\max_{k\in\left[  m\right]  ,\text{
}w_{\left\vert A\right\vert }^{r}\left(  k\right)  \neq0}\frac{w_{\left\vert
A\right\vert }^{r+p}\left(  k\right)  }{w_{\left\vert A\right\vert }%
^{r}\left(  k\right)  }%
\]
where $\left\vert A\right\vert =\left(  \left\vert a_{ij}\right\vert \right)
.$
\end{theorem}

The values $w_{A}^{r}$ can be used for lower bounds on $\sigma\left(
A\right)  $ as well.

\begin{theorem}
\label{th2} For every matrix $A$ and all $r\geq0$, $p\geq1,$
\[
\sigma^{2p}\left(  A\right)  \geq\frac{w_{A}^{r+p}}{w_{A}^{r}}%
\]
unless $\Sigma\left(  AA^{\ast}\right)  =0.$
\end{theorem}

On the other hand, for almost all matrices $A$ and $r$ large, Theorems
\ref{th1} and \ref{th2} are nearly optimal.

\begin{theorem}
\label{th3}For every $m\times n$ matrix $A$ and all $p\geq1,$
\[
\sigma^{2p}\left(  A\right)  =\lim_{r\rightarrow\infty}\frac{w_{A}^{r+p}%
}{w_{A}^{r}}=\lim_{r\rightarrow\infty}\max_{k\in\left[  m\right]  ,\text{
}w_{\left\vert A\right\vert }^{r}\left(  k\right)  \neq0}\frac{w_{A}%
^{r+p}\left(  k\right)  }{w_{A}^{r}\left(  k\right)  }%
\]
unless the eigenspace of $AA^{\ast}$ corresponding to $\sigma^{2}\left(
A\right)  $ is orthogonal to $\mathbf{j}_{m}$.
\end{theorem}

The following proposition sheds some light on Theorems \ref{th2} and \ref{th3}.

\begin{proposition}
\label{pro1}For every $m\times n$ matrix $A,$ the equality $w_{A}^{1}=0$ holds
if and only if $w_{A}^{r}=0$ holds for all $r\geq1.$ If $w_{A}^{1}=0,$ then
$\mathbf{j}_{m}$ is an eigenvector to $AA^{\ast}$ corresponding to $0;$
consequently all eigenvectors of $AA^{\ast}$ to nonzero eigenvalues are
orthogonal to $\mathbf{j}_{m}.$
\end{proposition}

Note also that, using Proposition \ref{pro3} below, Theorem \ref{th1} can be
extended to partitioned matrices. In particular, if $A$ is an $m\times n$
matrix partitioned into $pq$ blocks $A_{ij},$ $i\in\left[  p\right]  ,$
$j\in\left[  q\right]  ,$ then
\begin{equation}
\sigma^{2}\left(  A\right)  \leq\max_{i}%
{\textstyle\sum_{k=1}^{n}}
\sigma\left(  A_{ik}\right)
{\textstyle\sum_{k=1}^{m}}
\sigma\left(  A_{kj}\right)  \leq\max_{A_{ij}\neq0}%
{\textstyle\sum_{k=1}^{n}}
\sigma\left(  A_{ik}\right)
{\textstyle\sum_{k=1}^{m}}
\sigma\left(  A_{kj}\right)  . \label{partin}%
\end{equation}

\begin{proposition}
\label{pro3}Let the matrix $A$ be partitioned into $p\times q$ blocks
$A_{ij},$ $i\in\left[  p\right]  ,$ $j\in\left[  q\right]  .$ For all
$i\in\left[  p\right]  $ and $j\in\left[  q\right]  ,$ let $b_{ij}%
=\sigma\left(  A_{ij}\right)  .$ Then the matrix $B=\left(  b_{ij}\right)  $
satisfies $\sigma\left(  A\right)  \leq\sigma\left(  B\right)  .$
\end{proposition}

\section{Proofs}

\begin{proof}
[\textbf{Proof of Theorem \ref{th1}}]Since $\sigma\left(  A\right)  \leq
\sigma\left(  \left\vert A\right\vert \right)  ,$ to simplify the
presentation, we shall assume that $A$ is nonnegative. Likewise, dropping all
zero rows, we may assume that $A$ has no zero rows, that is to say, $w_{A}%
^{p}\left(  i\right)  >0$ for all $i\in\left[  m\right]  .$ Set $b_{ii}%
=w_{A}^{p}\left(  i\right)  $ for $i\in\left[  m\right]  $ and let $B$ be the
diagonal matrix with main diagonal $\left(  b_{11},\ldots,b_{mm}\right)  .$
Since $B^{-1}\left(  AA^{\ast}\right)  ^{r}B$ has the same spectrum as
$\left(  AA^{\ast}\right)  ^{k},$ the value $\sigma^{2r}\left(  A\right)  $ is
bounded from above by the maximum row sum of $B^{-1}\left(  AA^{\ast}\right)
^{r}B$ - say the sum of the $q$th row - and so,
\begin{align*}
\sigma^{2r}\left(  A\right)   &  \leq%
{\textstyle\sum\limits_{i\in\left[  m\right]  }}
w_{A}^{r}\left(  q,i\right)  \frac{w_{A}^{p}\left(  i\right)  }{w_{A}%
^{p}\left(  q\right)  }=\frac{1}{w_{A}^{p}\left(  q\right)  }%
{\textstyle\sum\limits_{i\in\left[  m\right]  }}
w_{A}^{r}\left(  q,i\right)
{\textstyle\sum\limits_{j\in\left[  m\right]  }}
w_{A}^{p}\left(  i,j\right) \\
&  =\frac{1}{w_{A}^{p}\left(  q\right)  }%
{\textstyle\sum\limits_{j\in\left[  m\right]  }}
{\textstyle\sum\limits_{i\in\left[  m\right]  }}
w_{A}^{r}\left(  q,i\right)  w_{A}^{p}\left(  i,j\right)  =\frac{1}{w_{A}%
^{p}\left(  q\right)  }%
{\textstyle\sum\limits_{j\in\left[  m\right]  }}
w_{A}^{p+r}\left(  q,j\right) \\
&  =\frac{w_{A}^{r+p}\left(  q\right)  }{w_{A}^{p}\left(  q\right)  }\leq
\max_{k\in\left[  m\right]  }\frac{w_{A}^{r+p}\left(  k\right)  }{w_{A}%
^{p}\left(  k\right)  },
\end{align*}
completing the proof.\bigskip
\end{proof}

\begin{proof}
[\textbf{Proof of inequalities (\ref{ourin})}]Theorem \ref{th1} with $r=0$ and
$p=1$ implies that%
\[
\sigma^{2}\left(  A\right)  \leq\max_{i\in\left[  m\right]  }w_{\left\vert
A\right\vert }^{1}\left(  i\right)  =\max_{i\in\left[  m\right]  }%
{\textstyle\sum\limits_{k\in\left[  m\right]  }}
{\textstyle\sum\limits_{j\in\left[  n\right]  }}
\left\vert a_{ij}\right\vert \left\vert a_{kj}\right\vert =\max_{i\in\left[
m\right]  }%
{\textstyle\sum\limits_{j\in\left[  n\right]  }}
\left\vert a_{ij}\right\vert
{\textstyle\sum\limits_{k\in\left[  m\right]  }}
\left\vert a_{kj}\right\vert =\max_{i\in\left[  m\right]  }%
{\textstyle\sum\limits_{a_{ij}\neq0}}
\left\vert a_{ij}\right\vert c_{j}.
\]
Suppose the maximum in the right hand side is attained for $i=k.$ Then,
\[%
{\textstyle\sum\limits_{j\in\left[  n\right]  }}
\left\vert a_{kj}\right\vert c_{j}=%
{\textstyle\sum\limits_{j\in\left[  n\right]  }}
\frac{\left\vert a_{kj}\right\vert }{r_{k}}r_{k}c_{j}\leq%
{\textstyle\sum\limits_{j\in\left[  n\right]  }}
\frac{\left\vert a_{kj}\right\vert }{r_{k}}\max_{a_{kj}\neq0}r_{k}c_{j}%
=\max_{a_{ij}\neq0}r_{i}c_{j},
\]
completing the proof in this case.
\end{proof}

\bigskip

In the proofs below we shall assume that $\sigma=\sigma_{1}\geq\cdots
\geq\sigma_{m}$ are the singular values of $A.$ Let $AA^{\ast}=VDV^{\ast}$ be
the unitary decomposition of $AA^{\ast};$ thus, the columns of $V$ are the
unit eigenvectors to $\sigma_{1}^{2},\ldots,\sigma_{m}^{2}$ and $D$ is the
diagonal matrix with $\sigma_{1}^{2},\ldots,\sigma_{m}^{2}$ along its main
diagonal. Writing $\Sigma\left(  B\right)  $ for the sum of the entries of a
matrix $B,$ note that for every $l\geq0,$
\[
w_{A}^{l}=\Sigma\left(  \left(  AA^{\ast}\right)  ^{l}\right)  =\Sigma\left(
VD^{l}V^{\ast}\right)  =%
{\textstyle\sum\limits_{i\in\left[  m\right]  }}
c_{i}\sigma_{i}^{2l},
\]
where $c_{i}=\left\vert \sum_{j\in\left[  m\right]  }v_{ji}\right\vert
^{2}\geq0$ is independent of $l.$\bigskip

\begin{proof}
[\textbf{Proof of Proposition \ref{pro1}}]In the notation above we see that
$w_{A}^{l}=0$ if and only if $c_{i}=0$ for every nozero $\sigma_{i},$ thus if
and only if $w_{A}^{1}=0.$

Note that $w_{A}^{1}=\Sigma\left(  AA^{\ast}\right)  =\left\langle AA^{\ast
}\mathbf{j}_{m},\mathbf{j}_{m}\right\rangle ;$ hence, if $w_{A}^{1}=0,$ then
$\mathbf{j}_{m}$ is an eigenvector of $AA^{\ast}$ to $0.$ Indeed, since
$AA^{\ast}$ is positive semidefinite, by the Rayleigh principle, $\left\langle
AA^{\ast}\mathbf{x},\mathbf{x}\right\rangle =0$ implies $AA^{\ast}%
\mathbf{x}=0.$ The proof is completed.
\end{proof}

\bigskip

\begin{proof}
[\textbf{Proof of Theorem \ref{th2}}]In the above notation we see that%
\[
\sigma^{2p}w_{A}^{2r}=%
{\textstyle\sum\limits_{i\in\left[  m\right]  }}
c_{i}\sigma^{2p}\sigma_{i}^{2r}\geq%
{\textstyle\sum\limits_{i\in\left[  m\right]  }}
c_{i}\sigma_{i}^{2p+2r}=w_{A}^{p+r}.
\]
The proof is completed by Proposition \ref{pro1}.\bigskip
\end{proof}

\begin{proof}
[\textbf{Proof of Theorem \ref{th3}}]Assume that there is an eigenvector of
$AA^{\ast}$ to $\sigma^{2}\left(  A\right)  $ that is not orthogonal to
$\mathbf{j}_{m}.$ Therefore, we may assume that $c_{1}>0.$ Hence,%
\[
\lim_{r\rightarrow\infty}\frac{%
{\textstyle\sum_{i\in\left[  m\right]  }}
c_{i}\sigma_{i}^{2p+2r}}{%
{\textstyle\sum_{i\in\left[  m\right]  }}
c_{i}\sigma_{i}^{2p}}=\sigma^{2p}\lim_{r\rightarrow\infty}\frac{%
{\textstyle\sum_{\sigma_{i}=\sigma_{1}}}
c_{i}}{%
{\textstyle\sum_{\sigma_{i}=\sigma_{1}}}
c_{i}}=\sigma^{2p},
\]
proving the first equality of the theorem.

For $k\in\left[  m\right]  $ and every $l\geq0,$ the value $w_{A}^{l}\left(
k\right)  $ is the $k$th row sum of the matrix $VD^{l}V^{\ast};$ hence
\[
w_{A}^{l}\left(  k\right)  =%
{\textstyle\sum\limits_{i\in\left[  m\right]  }}
{\textstyle\sum\limits_{j\in\left[  m\right]  }}
v_{ki}\sigma_{i}^{2l}\overline{v_{ji}}=%
{\textstyle\sum\limits_{i\in\left[  m\right]  }}
\sigma_{i}^{2l}v_{ki}%
{\textstyle\sum\limits_{j\in\left[  m\right]  }}
\overline{v_{ji}}=%
{\textstyle\sum\limits_{i\in\left[  m\right]  }}
b_{i}\sigma_{i}^{2l},
\]
where $b_{i}=v_{ki}%
{\textstyle\sum\limits_{j\in\left[  m\right]  }}
\overline{v_{ji}}$ is independent of $l.$ Writing $t$ for the largest number
such that $%
{\textstyle\sum_{\sigma_{i}=\sigma_{t}}}
b_{i}\neq0,$ we see that
\[
\lim_{r\rightarrow\infty}\frac{w_{A}^{r+p}\left(  k\right)  }{w_{A}^{r}\left(
k\right)  }=\sigma_{t}^{2p}\leq\sigma^{2p}.
\]
On the other hand, since%
\[
\max_{k\in\left[  m\right]  }\frac{w_{A}^{r+p}\left(  k\right)  }{w_{A}%
^{r}\left(  k\right)  }\geq\frac{%
{\textstyle\sum_{i\in\left[  m\right]  }}
w_{A}^{r+p}\left(  i\right)  }{%
{\textstyle\sum_{i\in\left[  m\right]  }}
w_{A}^{r+p}\left(  i\right)  },
\]
we obtain
\[
\liminf_{r\rightarrow\infty}\max_{k\in\left[  m\right]  }\frac{w_{A}%
^{r+p}\left(  k\right)  }{w_{A}^{r}\left(  k\right)  }\geq\sigma^{2p},
\]
completing the proof.\bigskip
\end{proof}

\begin{proof}
[\textbf{Proof of Proposition \ref{pro3}}]Let $A=\left(  a_{ij}\right)  $ be
an $m\times n$ matrix and $\left[  m\right]  =\cup_{i=1}^{p}P_{i}$ and
$\left[  n\right]  =\cup_{i=1}^{q}Q_{i}$ be the partitions of its index sets.
Select unit vectors $\mathbf{x}=\left(  x_{1},\ldots,x_{n}\right)  $ and
$\mathbf{y}=\left(  y_{1},\ldots,y_{m}\right)  $ such that $\sigma\left(
A\right)  =\left\langle A\mathbf{x},\mathbf{y}\right\rangle .$ Then we have
\begin{align*}
\sigma\left(  A\right)   &  =\left\langle A\mathbf{x},\mathbf{y}\right\rangle
=%
{\textstyle\sum\limits_{i\in\left[  m\right]  ,k\in\left[  n\right]  }}
a_{ik}x_{k}\overline{y_{i}}=%
{\textstyle\sum\limits_{r\in\left[  p\right]  }}
{\textstyle\sum\limits_{s\in\left[  q\right]  }}
{\textstyle\sum\limits_{i\in P_{r}}}
{\textstyle\sum\limits_{k\in Q_{s}}}
a_{ik}x_{k}\overline{y_{i}}\\
&  \leq%
{\textstyle\sum\limits_{r\in\left[  p\right]  }}
{\textstyle\sum\limits_{s\in\left[  q\right]  }}
\sigma\left(  A_{rs}\right)  \sqrt{%
{\textstyle\sum\limits_{i\in P_{r}}}
\left\vert x_{i}\right\vert ^{2}%
{\textstyle\sum\limits_{k\in Q_{s}}}
\left\vert y_{k}\right\vert ^{2}}\leq\sigma\left(  B\right)  ,
\end{align*}
completing the proof.\bigskip
\end{proof}

\textbf{Concluding remarks}

Theorem \ref{th1} and \ref{th2}\ extend Theorems 5 and 16 of \cite{Nik06},
that in turn generalize a number of results about the spectral radius of
graphs - see, e.g., the references of \cite{Nik06}.

Inequality (\ref{partin}) implies the essential result of the paper
\cite{DaBa05}; however, we admit that this paper triggered the present note.


\begin{thebibliography}{9}                                                                                                %


\bibitem {DaBa05}K.C. Das, R. Bapat, A sharp upper bound on the spectral
radius of weighted graphs, preprint available at
\emph{http://com2mac.postech.ac.kr/papers/2005/05-20.pdf}

\bibitem {HoJo88}R. Horn, C. Johnson, \emph{Matrix Analysis,} Cambridge
University Press, Cambridge, 1985. xiii+561 pp.

\bibitem {Nik06}V. Nikiforov, Walks and the spectral radius of graphs,
\emph{Linear Algebra Appl.} \textbf{418} (2006), 257-268.

\bibitem {Sch11}I. Schur, Bemerkungen zur Theorie der beschr\"{a}nkten
Bilinearformen mit unendlich vielen Ver\"{a}nderlischen, \emph{Journal f\"{u}r
Reine und Angew. Mathematik}, \textbf{140} (1911), 1--28.
\end{thebibliography}
\end{document}